\documentstyle[12pt,amstex,amssymb,amsthm]{article}

\newcommand{\C}{\mbox{\rm \,l\kern-0.52em C}}
\newcommand{\Ce}{\rm \,l\kern-0.35em C}
\newcommand{\R}{{\rm l}\!{\rm R}}

\newcommand{\N}{{\rm l}\!{\rm N}}

\newcommand{\Inf}{{\rm Inf}}
\newcommand{\Sup}{{\rm Sup}}

\newtheorem{theorem}{Theorem}[section]
\newtheorem{deflem}[theorem]{Definition and Lemma}
\newtheorem{prop}[theorem]{Proposition}

\newtheorem{cor}[theorem]{Corollary}

\newtheorem{example}[theorem]{Example}
\newtheorem{lemma}[theorem]{Lemma}
\newtheorem{remark}[theorem]{Remark}

\renewenvironment{proof}{{\bf Proof:}}{\mbox{}\hfill $\Box$}
\theoremstyle{definition}
\newtheorem{definition}[theorem]{Definition}
\title{New holomorphically closed subalgebras of $C^ *$-algebras of hyperbolic groups}

\author{Michael Puschnigg}
\date{}

\begin{document}
\maketitle
\section{Introduction}

The search for smooth subalgebras of $C^*$-algebras is an important problem in Noncommutative Geometry. Recall that a dense subalgebra of a Banach algebra is 
smooth (or closed under holomorphic functional calculus) if the spectra of its elements in both algebras agree. Typical examples are the algebra of smooth functions insisde the $C^*$-algebra of continuous functions on a smooth closed manifold, or the Banachalgebra of trace-class operators inside the $C^*$-algebra of all compact operators on a Hilbert space.

Certain invariants of Banach algebras like $K$-groups or local cyclic cohomology groups  do not change under passage to a smooth subalgebra. The construction of a particularly small and well behaved smooth subalgebra is therefore often an important first step towards the calculation of such invariants. 

In this paper we construct new classes of smooth subalgebras of the reduced group 
$C^*$-algebra $C^*_r(\Gamma)$ of a word-hyperbolic group $\Gamma$. In particular, we find a smooth unconditional subalgebra ${\cal A}_\omega\Gamma$ which possesses the following property: every trace on the group ring $\C\Gamma$, which is supported in finitely many conjugacy classes of $\Gamma$, extends to a bounded trace on ${\cal A}_\omega\Gamma$. This result has two consequences. 

The first concerns the cyclic homology of group $C^*$-algebras 
and of unconditional completions of group rings. Recall that a Banach- or Fr\'echetalgebra over the group ring is unconditional \cite{La}, if the norm of an element $\underset{g}{\sum}a_g\,u_g\,\in\,\C\Gamma$ 
only depends on the absolute values $(\vert a_g\vert)_{g\in\Gamma}$ of its Fourier coefficients. 

It is well known \cite{Ni} that cyclic complexes of group rings of discrete groups decompose as  direct sums of complexes which are labeled by the conjugacy classes. A similar (topological) decomposition exists for the analytic cyclic complex of the group Banach algebra $\ell^1(\Gamma)$, but not for other completions of the group ring. 

Here we establish a weak form of 
homogeneous decomposition for the analytic cyclic complex of any unconditional completion of the group ring of a word-hyperbolic group. The contribution of a fixed conjugacy class turns out to be independent of the unconditional completion under study. For the conjugacy class of the unit this result is due to Connes and Moscovici \cite{CM} and plays a crucial role in their proof of the Novikov higher signature conjecture for word-hyperbolic groups. For us the homogeneous decomposition is the first basic step towards a calculation of the bivariant local cyclic cohomology of unconditional Banachalgebras and the reduced $C^*$-algebra of a hyperbolic group. This will be the subject of another paper.

Another application of our result concerns a problem of John Lott \cite{Lo1}, \cite{Lo2}. As Thomas Schick pointed out to me, the results of John Lott on delocalized $L^2$-invariants of smooth closed manifolds with hyperbolic fundamental group depend on the existence of smooth subalgebras of the type we constructed here. So the rather restrictive conditions imposed in \cite{Lo2} can be removed and Lott's results hold as originally stated in \cite{Lo1}. 

Our construction of the smooth subalgebras ${\cal A}_\omega\Gamma$ of $C^*_r(\Gamma)$ is based on two things. 

The first is a particular type of topological tensor product for unconditional Fr\'echet spaces.
It is easily seen that the class of unconditional Banach- or Fr\'echetalgebras over group rings is not closed under (projective) tensor products. This is the source of trouble in the calculation of cyclic homology for such an algebra because the cyclic bicomplex is a direct sum of projective tensor powers of the algebra. 

We show that there exists a largest unconditional cross-norm on the algebraic tensor product of unconditionally normed spaces and derive various properties of this unconditional tensor product. 

The second invention is a linear operator 
$\Delta:\C\Gamma\to\C\Gamma\otimes\C\Gamma$, which we call the canonical quasiderivation (despite the fact that it depends on choices). If one fixes a finite symmetric set $S$ of generators of $\Gamma$ and chooses a word $(s_1,...,s_n)$ of minimal length (in the alphabet $S$) for $g\in\Gamma$, then
$$
\begin{array}{ccc}
\Delta(u_g) & = & \underset{i=0}{\overset{n}{\sum}}\,u_{s_1...s_i}\otimes u_{s_{i+1}...s_n}
\end{array}
$$
It turns out that $\Delta$ behaves in some sense like a derivation, provided that $\Gamma$ is word-hyperbolic. If $\Delta$ is viewed as a densely defined unbounded operator on an unconditional Banach algebra ${\cal A}\Gamma$ over $\C\Gamma$ with values in the unconditional tensor square ${\cal A}\Gamma^{\otimes^2_{uc}}$, then the graph norm   is a differential Banach norm in the sense of Blackadar and Cuntz \cite{BC}. In particular, the closure of the group ring in the graph norm is a dense and holomorphically closed subalgebra of ${\cal A}\Gamma$. 
It should be noted that the proof depends strongly on hyperbolicity and that there is no indication that a similar result holds for other classes of groups.
Applying a suitable iteration of this procedure to the Jolissaint algebra \cite{Ha}, \cite{Jol} of square summable functions on $\Gamma$ leads to the desired algebras ${\cal A}_\omega\Gamma$.

A crucial property of word-hyperbolic groups (which was established by Gromov and actually characterizes this class of groups) is the existence of a linear isoperimetric inequality \cite{Gr}, 2.3. Closely related to this is the existence of solutions of linear growth to the conjugacy problem in word-hyperbolic groups \cite{Gr}, 7.4. This means that for two conjugate elements $g,g'$ in a word-hyperbolic group one may find an element $u$ satisfying $g'=ugu^{-1}$ and such that its word length is bounded from above by a fixed linear function of the word lengths of $g$ and $g'$. Making use of this result of Gromov one sees easily that the restriction of functions on a word-hyperbolic group $\Gamma$ to a fixed conjugacy class $\langle x\rangle$ gives rise to a bounded operator 
$$
\begin{array}{cccc}
 res_{\langle x\rangle}: & {\cal A}_\omega\Gamma & \to & \ell^1(\Gamma).
\end{array}
$$
This leads immediately to the consequences described before.

The results of the present paper generalize those obtained by different methods for free groups in \cite{Pu5}. Nevertheless the solution of the special case was crucial for finding the methods employed here.

I thank Thomas Schick heartly for discussions and for bringing John Lott's work to my attention.

\section{Unconditional completions}

\subsection{Unconditional seminorms}

Let $X$ be a set and denote by $\C X$ the complex vector space with basis $X$. It can be identified with the space of complex measures on $X$ with finite support. Any element can be uniquely written as 
$$
a=\sum a_xu_x\,\in\C X
\eqno(2.1)
$$ 
Its absolute value is defined as 
$$
\vert a\vert=\sum\vert a_x\vert u_x\,\in\C X
\eqno(2.2)
$$
We introduce a partial order on {\bf radial elements} $a=\vert a\vert,\,b=\vert b\vert\in\C X$ by
$$ 
\begin{array}{cccc}
a\leq b & \Leftrightarrow & a_x\leq b_x & \forall \,x\in X
\end{array}
\eqno(2.3)
$$

\begin{definition}
A linear map $\varphi:\C X\to\C Y$ is {\bf radial} if 
$$
\begin{array}{ccc}
\vert\varphi(a)\vert & \leq & \varphi(\vert a\vert) \\
\end{array}
\eqno(2.4)
$$
for all $a\in\C X$.
\end{definition}
In particular, a radial map sends radial elements to radial elements and preserves the partial order on the set of radial elements.
\begin{lemma}
Radial maps are stable under composition.
\end{lemma}

\begin{proof}
Let $f:\C X\to\C Y$ and $g:\C Y\to\C Z$ be radial maps. Because radial linear maps preserve the partial order on radial elements, one finds for $a\in\C X$ that\\
$\vert (g\circ f)(a)\vert\leq g(\vert f(a)\vert)\leq g(f(\vert a\vert))$ as was to be shown.
\end{proof}

The following notion is due to Bost and Lafforgue \cite{La}.
\begin{definition}
A seminorm on $\C X$ is called {\bf unconditional} if 
$$
\begin{array}{cccc}
\vert a\vert\,\leq\,\vert b\vert & \Rightarrow & \parallel a\parallel\,\leq\,\parallel b\parallel, & \forall \,a,b\in\C X
\end{array}
\eqno(2.5)
$$
 A Fr\'echet space is called an {\bf unconditional completion} of $\C X$ if it contains $\C X$ as a dense subspace and can be defined by a countable family of seminorms which are unconditional on $\C X$.
\end{definition}

Typical examples of unconditional seminorms are (weighted) $\ell^p$-norms, $1\leq p\leq\infty$. On the contrary, if $X=\Gamma$ is a discrete group with more than one element, then the $C^*$-norms of the maximal or reduced group $C^*$-algebra are not unconditional.

\begin{deflem}
Let $\parallel-\parallel$ be a seminorm on $\C X$. Then
$$
\begin{array}{cccc}
\parallel\alpha\parallel^{+} & := & \underset{\vert\alpha\vert\leq\vert\beta\vert}{\Inf}\,
\parallel\vert\beta\vert\parallel, & \alpha\in\C X
\end{array}
\eqno(2.6)
$$
defines the {\bf associated unconditional seminorm} on $\C X$. It is the largest unconditional seminorm which is majorized on radial elements by the originally given seminorm.
\end{deflem}

This construction is not very useful because in most cases there is no relation between a norm 
and its associated unconditional norm. A notable exception is the unconditional cross norm 
on tensor products which we will  discuss in a moment. 

Most norm estimates in this section will be based on 

\begin{lemma}
Let $X,Y$ be sets and let $\parallel-\parallel_X,\,\parallel-\parallel_Y$ be seminorms on $\C X,\,\C Y$, 
respectively. Let $\varphi:(\C X,\,\parallel-\parallel_X)\to(\C Y,\,\parallel-\parallel_Y)$ be a bounded and radial linear map. Then $\varphi$ is also bounded with respect to the corresponding unconditional seminorms $\parallel-\parallel^{+}_X,\,\parallel-\parallel^{+}_Y$
and 
$$
\begin{array}{ccc}
\parallel\varphi\parallel^{+} & \leq & \parallel\varphi\parallel \\
\end{array}
\eqno(2.7)
$$
\end{lemma}

\begin{proof}
Let $\alpha,\beta\in\C X$ be such that $\vert\alpha\vert\leq\vert\beta\vert$.
Then $\vert\varphi(\alpha)\vert\leq\varphi(\vert\alpha\vert)\leq\varphi(\vert\beta\vert)$ because $\varphi$ is radial and radial maps preserve the partial order on radial elements. Therefore
$$
\parallel\varphi(\alpha)\parallel^{+}\,=\,
\underset{\vert\varphi(\alpha)\vert\leq\vert\gamma\vert}{\Inf}\,
\parallel\vert\gamma\vert\parallel\,\leq\,
\underset{\vert\alpha\vert\leq\vert\beta\vert}{\Inf}\,
\parallel\varphi(\vert\beta\vert)\parallel
\leq\,\parallel\varphi\parallel
\underset{\vert\alpha\vert\leq\vert\beta\vert}{\Inf}\,
\parallel\vert\beta\vert\parallel
\,=\,\parallel\varphi\parallel\parallel\alpha\parallel^{+}
$$
\end{proof}

\subsection{Unconditional completions of tensor products} 

A seminorm on the algebraic tensor product $V\otimes W$ of seminormed vector spaces $V,W$ is a {\bf cross-seminorm} if 
$$
\begin{array}{cc}
\parallel v\otimes w\parallel\,=\,\parallel v\parallel_V\cdot\parallel w\parallel_W, & \forall\,v\in V,\,\forall\,w\in W. \\
\end{array}
\eqno(2.8)
$$
 A cross-seminorm is by no means determined by the seminorms on the individual factors, but there exists a largest  cross-seminorm, called the projective cross seminorm. The completion of the algebraic tensor product of Fr\'echet spaces with respect to the projective cross-seminorms is denoted by $V\otimes_\pi W$. It is characterized by an obvious universal property. 

It is not true that the class of unconditional seminorms is stable under passage to projective tensor products.

\begin{example}
Let $X=\{x,y\}$ be a set with two elements and consider the projective cross-norm on $\C(X\times X)\,\simeq\,\C X\otimes\C X$ associated to the $\ell^2$ norm on $\C X$. The latter is obviously unconditional.
Put 
$$a\,=\,u_{(x,x)}\,+\,2u_{(x,y)}\,+\,2u_{(y,x)}\,\in\,\C(X\times X)$$ and 
$$b\,=\,u_{(x,x)}\,+\,2u_{(x,y)}\,+\,2u_{(y,x)}\,+\,u_{(y,y)}\,\in\,\C(X\times X).$$
Then 
$$0\,\leq\,a\,\,\leq\,b\,\in\,\C(X\times X)$$ but  
$$
\parallel a\parallel_{\ell^2\otimes_\pi\ell^2}=
\sqrt{17}\,>\,4\,=\,\parallel b\parallel_{\ell^2\otimes_\pi\ell^2}
$$
\end{example}
In fact, if $\cal H$ is a separable Hilbert space (in our case $\ell^2(X)$), then ${\cal H}\otimes_\pi{\cal H}\,\simeq\,{\cal H}^*\otimes_\pi{\cal H}\,\subset\,{\cal L(H)}$ coincides with the space $\ell^1({\cal H})$ of trace class operators on $\cal H$. The norm of a trace class operator is given by the sum of its singular values. In the case of the elements constructed above the corresponding operators are given by $2\times 2$-matrices and their singular values can be calculated explicitely.

This fact is a serious obstacle if one wants, for example, to calculate the cyclic homology of an unconditional Fr\'echet algebra associated to a discrete grupoid. The cyclic chain complex consists of a direct sum of all tensor powers of the given algebra and it becomes very difficult to establish the necessary estimates if the property of being unconditional is lost. A striking example is given by the Jolissaint-algebra \cite{Jol}, which is a dense and holomorphically closed subalgebra of the reduced $C^*$-algebra of a word-hyperbolic group. It is defined by a family of weighted $\ell^2$-norms on the group ring, so that under passage to tensor powers one encounters exactly the phenomenon described in the previous example. 

\begin{prop}
Let $X,Y$ be sets and let $\parallel-\parallel_X, \parallel-\parallel_Y$ be unconditional seminorms on $\C X,\,\C Y$, respectively. Then the unconditional seminorm $\parallel-\parallel_{uc}$ \\ associated (see 2.4) to the projective cross-norm on $\C X\otimes\C Y\,\simeq\,\C(X\times Y)$ is given by
$$
\begin{array}{cccc}
\parallel\alpha\parallel_{uc} & = & \underset{I}{\Inf}\,\,\underset{i\in I}{\sum}\parallel\alpha_i'\parallel_X\parallel\alpha_i''\parallel_Y, & 
\forall\alpha\in\C(X\times Y) \\
\end{array}
\eqno(2.9)
$$
where the infimum is taken over all finite families $\alpha_i'\in\C X,\,\alpha_i''\in\C Y,\,i\in I,$ satisfying 
$$
\begin{array}{ccc}
\vert\alpha\vert & \leq & \underset{i\in I}{\sum}\,\,\vert\alpha_i'\vert\otimes\vert\alpha_i''\vert. \\
\end{array}
\eqno(2.10)
$$
It is the largest unconditional cross-seminorm on $\C X\otimes\C Y\,\simeq\,\C(X\times Y)$.

If $E,F$ are unconditional completions of $\C X,\,\C Y$, respectively, then the completion of the algebraic tensor product $E\otimes F$ with respect to the 
unconditional cross-seminorms (2.9) is denoted by $E\otimes_{uc}F$.
\end{prop}

\begin{proof}
Let $\alpha,\,\beta\in\C(X\times Y)$ and let $\alpha_i',\beta_j'\in\C X,\,\alpha_i'',\beta_j''\in\C Y,\,i\in I,\,j\in J$ be finite families of elements such that $\vert\alpha\vert\leq\underset{i\in I}{\sum}\,\,\vert\alpha_i'\vert\otimes\vert\alpha_i''\vert$ and 
$\vert\beta\vert\leq\underset{j\in J}{\sum}\,\,\vert\beta_j'\vert\otimes\vert\beta_j''\vert$, respectively. Then 
$$
\vert\alpha+\beta\vert\,\leq\,\vert\alpha\vert+\vert\beta\vert\,\leq\,\underset{i\in I}{\sum}\,\,\vert\alpha_i'\vert\otimes\vert\alpha_i''\vert
\,+\,\underset{j\in J}{\sum}\,\,\vert\beta_j'\vert\otimes\vert\beta_j''\vert
$$
so that by definition 
$$
\parallel\alpha+\beta\parallel\,\leq\,\underset{I,J}{\Inf}\,\,\left(
\underset{i\in I}{\sum}\parallel\alpha_i'\parallel_X\parallel\alpha_i''\parallel_Y\,+\,
\underset{j\in J}{\sum}\parallel\beta_j'\parallel_X\parallel\beta_j''\parallel_Y\right)\,
\,=\,\parallel\alpha\parallel\,+\,\parallel\beta\parallel
$$
Thus the functional (2.9) is subadditive. Its homogeneity is clear. We show next that it is a cross-seminorm.
 Let $\xi\in\C X,\,\eta\in\C Y$. The identity $\vert\xi\otimes\eta\vert=\vert\xi\vert\otimes\vert\eta\vert$ shows already that $\parallel\xi\otimes\eta\parallel\leq\parallel\xi\parallel\parallel\eta\parallel$. Put now $X'=Supp(\xi)=\{x_1,\ldots,x_n\}\subset X,\,Y'=Supp(\eta)
 =\{y_1,\ldots,y_m\}\subset Y$. By the theorem of Hahn-Banach there exist linear functionals 
 $l:\R X'\to\R,\,l':\R Y'\to\R$ of norm one such that \\ $l(\vert\xi\vert)=\parallel\xi\parallel,\,l'(\vert\eta\vert)=\parallel\eta\parallel$. If $\xi=\lambda_1x_1+\ldots+\lambda_nx_n$ we put $\xi_i\,=\,\underset{j\neq i}{\sum}\,\lambda_j x_j$ for $1\leq i\leq n$. Then 
$l(\vert\xi_i\vert)\leq\vert l(\vert\xi_i\vert)\vert\leq\parallel\vert\xi_i\vert\parallel\,
\leq\,\parallel\vert\xi\vert\parallel\,=\,\parallel\xi\parallel\,=\,l(\vert\xi\vert)$ because the seminorm is unconditional and thus $l(x_i)\,=\,\vert\lambda_i\vert^{-1}l(\vert\xi\vert-\vert\xi_i\vert)\geq 0$. Similarly 
$l'(y)\geq 0$ for $y\in Y'$ so that $\vert\alpha\vert\leq\vert\beta\vert\Rightarrow
l\otimes l'(\vert\alpha\vert)\leq l\otimes l'(\vert\beta\vert)$ for $\alpha,\,\beta\in\C(X'\times Y')$.
Let finally be $m:\C X\to\C X'$ and $m':\C Y\to\C Y'$ be the operators given by multiplication with the characteristic functions of $X'$ and $Y'$, respectively. They are of norm one because the given seminorms are unconditional. Let now $\alpha_i'\in\C X,\,\alpha_i''\in\C Y,\,i\in I$ be a finite family such that $\vert\xi\otimes\eta\vert \leq\underset{i\in I}{\sum}\,\,\vert\alpha_i'\vert\otimes\vert\alpha_i''\vert$.
Then $\vert\xi\otimes\eta\vert \leq\underset{i\in I}{\sum}\,\,\vert m(\alpha_i')\vert\otimes\vert m'(\alpha_i'')\vert\in\C(X'\times Y')$ and 
$\parallel\xi\parallel\parallel\eta\parallel=l(\vert\xi\vert)l'(\vert\eta\vert)$
$\leq(l\otimes l')(\underset{i\in I}{\sum}\,\,\vert m(\alpha_i')\vert\otimes\vert m'(\alpha_i'')\vert)$
$=\underset{i\in I}{\sum}\,\,l(\vert m(\alpha_i')\vert)l'(\vert m'(\alpha_i'')\vert)$
$\leq\underset{i\in I}{\sum}\,\,\parallel\vert m(\alpha_i')\vert\parallel\parallel\vert m'(\alpha_i'')\vert\parallel$
$\leq\underset{i\in I}{\sum}\,\,\parallel\alpha_i'\parallel\parallel\alpha_i''\parallel$
which shows finally that $\parallel\xi\otimes\eta\parallel=\parallel\xi\parallel\parallel\eta\parallel$.

Let now $\parallel-\parallel'$ be any unconditional norm on $\C X\otimes\C Y$ which is majorized by the projective cross-norm on radial elements. Let $\beta\in\C(X\times Y)$ and let $\beta_j'\in\C X,\,\beta_j''\in\C Y,\,j\in J$ be a finite family such that $\vert\beta\vert \leq  \underset{j\in J}{\sum}\,\,\vert\beta_j'\vert\otimes\vert\beta_j''\vert$. Then 
$$
\parallel\beta\parallel'\,=\,\parallel\vert\beta\vert\parallel'\,\,
\leq\,\,
\parallel\underset{j\in J}{\sum}\,\,\vert\beta_j'\vert\otimes\vert\beta_j''\vert\parallel'
\,\,\leq\,\,\parallel\,\underset{j\in J}{\sum}\,\vert\beta_j'\vert\otimes\vert\beta_j''\vert\parallel_\pi\,\,\leq
$$
$$
\leq\,\underset{j\in J}{\sum}\,\,\parallel\vert\beta_j'\vert\parallel_X\cdot\parallel\vert\beta_j''\vert\parallel_Y\,\,
=\,\,\underset{j\in J}{\sum}\,\,\parallel\beta_j'\parallel_X\cdot\parallel\beta_j''\parallel_Y
$$
which shows that (2.9) is the largest unconditional cross norm which is majorized 
by the projective cross-norm on radial elements. It coincides therefore with the unconditional seminorm $\parallel-\parallel_{uc}$ associated (see 2.4) to the projective cross-norm on $\C X\otimes\C Y$.
\end{proof}

\subsection{Unconditional completions of algebras}

\begin{definition}
\begin{itemize}
\item[a)]
Let $X$ be a set. Then a complex algebra structure on $\C X$ is called {\bf radial}
if the multiplication $m:\C X\otimes\C X\,=\,\C(X\times X)\to \C X$ 
is radial in the sense of 2.1. 
\item[b)]
A Banach(Fr\'echet)algebra obtained from a radial algebra $\C X$ by completion 
with respect to a family of unconditional seminorms is called an \\ {\bf unconditional 
Banach(Fr\'echet) algebra over $\C X$}.
\end{itemize}
\end{definition}
A class of examples of radial algebras is provided by
\begin{lemma}
Let $\cal C$ be a small category and let $Mor({\cal C})$ be the set of its morphisms.
Then the convolution algebra $\C Mor({\cal C})$ is radial. 
\end{lemma}
\begin{proof}
This is immediate from
$$
\vert m(\alpha)\vert(g)=\vert\underset{g'g''=g}{\sum}\alpha(g',g'')\vert
\leq\underset{g'g''=g}{\sum}\vert\alpha(g',g'')\vert=m(\vert\alpha\vert)(g),
$$
which holds for all $\alpha\in\C(Mor({\cal C})\times Mor({\cal C}))$.
\end{proof}
Examples of unconditional topological algebras over the group ring $\C\Gamma$ 
of a discrete group are the group Banach algebra $\ell^1(\Gamma)$, the Jolissaint algebra (6.2),\cite{Ha}, \cite{Jol}, in the case that $\Gamma$ is hyperbolic, and the completions of the group ring with respect to unconditional norms associated to submultiplicative seminorms via 2.4.

\begin{lemma}
Let ${\cal A}$ be an unconditional Fr\'echet algebra over the radial algebra $\C X$. Then the multiplication on ${\cal A}$ extends to 
a bounded operator
$$
\begin{array}{cccc}
m: & {\cal A}\otimes_{uc}{\cal A} & \longrightarrow & {\cal A} \\
\end{array}
\eqno(2.11)
$$
More generally, if $V_Y$ is an unconditional completion of  $\C_Y$ for some set $Y$, then the multiplication 
on $A$ gives rise to a bounded linear map
$$
\begin{array}{cccc}
m\otimes id: & {\cal A}\otimes_{uc}{\cal A}\otimes_{uc} V_Y & \longrightarrow & 
{\cal A}\otimes_{uc}V_Y \\
\end{array}
\eqno(2.12)
$$
The norms of these maps are majorized by their norms on the corresponding projective tensor products.
\end{lemma}
\begin{proof}
By definition the multiplication on $\C X$ is radial. Therefore the assertion follows from
 2.5 and 2.7. 
\end{proof}

\begin{lemma}
Let $\C X$ be a radial algebra and let $T\subset\C X$ be a subset generating $\C X$ as complex algebra. There exists a largest submultiplicative unconditional seminorm 
$\parallel-\parallel_T$ on $\C X$ such that $\parallel T\parallel_T\leq 1$.

For $\alpha\in\C X$ it is given by 
$$
\parallel\alpha\parallel_T  =  \underset{I}{\Inf}\,\,\underset{i\in I}{\sum}\,\vert\lambda_i\vert
$$
where the infimum is taken over all finite families $t_{i,1},\ldots , t_{i,n_i}\in T,\,\lambda_i\in\C,\,i\in I,$ satisfying 
$$
\begin{array}{ccc}
\vert\alpha\vert & \leq & \underset{i\in I}{\sum} \vert\lambda_i\vert\vert t_{i,1}\vert\ldots\vert t_{i,n_i}\vert.
\end{array}
\eqno(2.13)
$$
In particular, if $T$ is stable under $a\mapsto\vert a\vert$, then $\parallel-\parallel_T$ coincides with the largest unconditional seminorm such that $\parallel T^n\parallel_T\,\leq\, 1$ 
for all $n\in\N$.
\end{lemma}

\section{The canonical quasiderivation}

\begin{definition}
Let $\Gamma$ be a finitely generated group and let $S=S^ {-1}$ be a finite symmetric set of generators.
Let $\sigma:\Gamma\to\underset{n}{\bigcup}S^n$ be a map which associates to every element $g\in\Gamma$ a word of minimal length (in the alphabet $S$) representing $g$. (Such a map will be called a minimal presentation of the elements of $\Gamma$.)

The linear operator 
$$
\begin{array}{cccc}
\Delta_{(S,\sigma)}: & \C\Gamma & \to & \C\Gamma\otimes\C\Gamma   \\
 & & & \\
 & u_g & \mapsto & \underset{i=0}{\overset{n}{\sum}}\,u_{s_1\ldots s_{i}}\otimes u_{s_{i+1}\ldots s_n}
\end{array}
\eqno(3.1)
$$
where
$$
\sigma(g)=(s_1,\ldots,s_n)
$$
is called the {\bf canonical quasiderivation} of $\C\Gamma$ associated to $(S,\sigma)$. 
\end{definition}

Let ${\cal G}(\Gamma,S)$ be the Cayley graph of $\Gamma$ with respect to $S$. The set of vertices of ${\cal G}(\Gamma,S)$ coincides with $\Gamma$ and two vertices $g,\,g'\in\Gamma$ are adjacent iff $g^ {-1}g'\in S$. The maximal metric which identifies edges with the unit interval turns the Cayley graph into a geodesic metric space. The left multiplication induces
 a simplicial and isometric action of $\Gamma$ on ${\cal G}(\Gamma,S)$. 
Every word of minimal length for an element $g\in G$ determines a geodesic path from 
$h\in{\cal G}(\Gamma,S)^ 0$ to $hg\in{\cal G}(\Gamma,S)^ 0$. Thus a minimal presentation of the elements of $\Gamma$ can be interpreted as a translation invariant choice of geodesic segment between any two vertices of ${\cal G}(\Gamma,S)$.

Recall that $\Gamma$ is {\bf word-hyperbolic} in the sense of Gromov \cite{Gr} if for some 
(and thus every) choice of a finite symmetric set $S$ of generators the Cayley graph ${\cal G}(\Gamma,S)$ is a hyperbolic metric space in the following sense:

\begin{definition} \cite{Gr}, 2.1.A, 6.3.C.
A geodesic metric space $(X,d)$ is {\bf hyperbolic} if there exists a constant $\delta\geq 0$ such that each edge of a geodesic triangle in $X$ is contained in a $\delta$-neighbourhood of the union of the two other edges. 
\end{definition}
Thus, if $a,\,b,\,c$ are the edges of a geodesic triangle, then 
$$
\begin{array}{ccc}
 x\in a & \Rightarrow & d_X(x,\, b\cup c)\,\leq\,\delta \\
\end{array}
\eqno(3.2)
$$
In the sequel we suppose (without loss of generality) that $\delta>0$ is an integer.

\begin{theorem}
Let $\Gamma$ be a word-hyperbolic group \cite{Gr}. Choose a finite symmetric set of generators $S$, denote by $\sigma$ be an associated minimal presentation of the elements of $\Gamma$, and let $\Delta_{(S,\sigma)}:\C\Gamma\to
\C\Gamma\otimes\C\Gamma$ be the associated quasiderivation. Let ${\cal A}\Gamma$ be an unconditional Banach algebra over $\C\Gamma$ (see 2.8) with defining seminorm $\parallel-\parallel$.
Then there exists a constant $C_0$ (depending on all chosen data) such that
$$
\begin{array}{ccc}
\parallel\Delta_{(S,\sigma)}(\alpha\beta)\parallel_{uc} & \leq & 
C_0\left( \parallel\Delta_{(S,\sigma)}(\alpha)\parallel_{uc}\cdot\parallel\beta\parallel\,+\,
\parallel\alpha\parallel\cdot\parallel\Delta_{(S,\sigma)}(\beta)\parallel_{uc}\right)
\end{array}
\eqno(3.3)
$$
for all $\alpha,\,\beta\in\C\Gamma$.
\end{theorem}

It should be noted that it is essential to take the unconditional instead of the projective cross norm.

\begin{proof}
For $g,\,g'\in\Gamma$ consider the geodesic triangle in $\cal G(\Gamma,S)$ with edges $[1,g]=\sigma(g),\,[g,gg']=g\sigma(g')$ and $[1,gg']=\sigma(gg')$.
For every vertex $h\in [1,gg']$ one can either find a vertex $h_0\in[1,g]$ with $d(h,h_0)<\delta$ or a vertex $h_1\in[1,g']$ with $d(h,gh_1)<\delta$. Thus 
$$
\Delta_{(S,\sigma)}(u_gu_{g'})\,=\,\underset{h\in[1,gg']}{\sum}\,u_h\otimes u_{h^{-1}gg'}
$$
$$
\leq\,\underset{l_S(t)\leq\delta}{\sum}\left(\underset{h_0\in[1,g]}{\sum}\,u_{h_0t}\otimes 
u_{t^{-1}(h_0^{-1}g)} u_{g'}\,
+\,\underset{h_1\in[1,g']}{\sum}\,u_g u_{h_1t}\otimes 
u_{t^{-1}(h_1^{-1}g')}\right)
\eqno(3.4)
$$
Let now $\alpha=\sum a_gu_g,\,\beta=\sum b_{g'}u_{g'}$ be elements of $\C\Gamma$.
Then we deduce from (3.4)
$$
\vert\Delta_{(S,\sigma)}(\alpha\beta)\vert\,\leq\,
\underset{g,g'}{\sum}\vert a_g\vert\vert b_{g'}\vert\underset{h\in[1,gg']}{\sum}\,u_h\otimes u_{h^{-1}gg'}
$$
$$
\leq\,\underset{g,g'}{\sum}\vert a_g\vert\vert b_{g'}\vert\underset{l_S(t)\leq\delta}{\sum}\left(
\underset{h_0\in[1,g]}{\sum}\,u_{h_0t}\otimes 
u_{t^{-1}(h_0^{-1}g)} u_{g'}\,
+\,\underset{h_1\in[1,g']}{\sum}\,u_g u_{h_1t}\otimes 
u_{t^{-1}(h_1^{-1}g')}\right)
$$
$$
=\underset{l_S(t)\leq\delta}{\sum}\,i_t\left(\vert\Delta_{(S,\sigma)}(\alpha)\vert\vert\beta\vert\,+\,
\vert\alpha\vert\vert\Delta_{(S,\sigma)}(\beta)\vert
\right)
$$
$$
=\underset{l_S(t)\leq\delta}{\sum}\,\left(i_t(\vert\Delta_{(S,\sigma)}(\alpha)\vert)\vert\beta\vert\,+\,\vert\alpha\vert i_t(\vert\Delta_{(S,\sigma)}(\beta)\vert)
\right)
\eqno(3.5)
$$
where $i_t$ is the $\C\Gamma$-bimodule endomorphism of $\C\Gamma\otimes\C\Gamma$ given by 
$$
\begin{array}{cccccc}
i_t & = & m_r(u_{t^{-1}})\otimes m_l(u_t): & u_g\otimes u_{g'} & \mapsto &
u_{gt^{-1}}\otimes u_{tg'}
\end{array}
\eqno(3.6)
$$
Note that it is a radial map so that, according to 2.5 and 2.7 
$$
\parallel i_t\parallel_{uc}\,\leq\,\parallel i_t\parallel_{\pi}\,=\,\parallel m_r(u_t)\otimes_\pi m_l(u_{t^{-1}})\parallel_{\pi}\,=\,\parallel u_t\parallel\cdot\parallel u_{t^{-1}}\parallel
\eqno(3.7)
$$
Therefore  
$$
\parallel\Delta_{(S,\sigma)}(\alpha\beta)\parallel_{uc}\,\leq\,
\parallel\underset{l_S(t)\leq\delta}{\sum}\,i_t\left(\vert\Delta_{(S,\sigma)}(\alpha)\vert\vert\beta\vert\,+\,
\vert\alpha\vert\vert\Delta_{(S,\sigma)}(\beta)\vert
\right)\parallel_{uc}
$$
$$
\leq\left(\underset{l_S(t)\leq\delta}{\sum}\parallel u_t\parallel\cdot\parallel u_{t^{-1}}\parallel\right)
\left(\parallel\Delta_{(S,\sigma)}(\alpha)\parallel_{uc}\cdot\parallel\beta\parallel\,+\,
\parallel\alpha\parallel\cdot\parallel\Delta_{(S,\sigma)}(\beta)\parallel_{uc}\right)
$$
which proves the theorem with the constant 
$$
C_0\,=\,\underset{l_S(t)\leq\delta}{\sum}\parallel u_t\parallel\cdot\parallel u_{t^{-1}}\parallel.
\eqno(3.8)
$$ 
\end{proof}

\begin{definition}
Let $(\Gamma,S)$ be a word-hyperbolic group and let $\sigma$ be a minimal presentation of its elements in terms of the alphabet $S$. Then for an unconditional Banach algebra ${\cal A}\Gamma$ over $\C\Gamma$ we denote by 
${\cal A}^1\Gamma$ the closure of $\C\Gamma$ in the graph-norm of the canonical quasiderivation $\Delta_{(S,\sigma)}$, viewed as unbounded linear operator from ${\cal A}\Gamma$ to ${\cal A}\Gamma\otimes_{uc}{\cal A}\Gamma$.
\end{definition}

The previous theorem shows that the graph norm of the canonical quasiderivation is 
a differential norm of "logarithmic order $C_0$" in the sense of Blackadar and Cuntz \cite{BC}, provided that the underlying group is hyperbolic. Consequently

\begin{theorem}
Let $\Gamma$ be a word-hyperbolic group and let ${\cal A}\Gamma$ be an unconditional Banach algebra over $\C\Gamma$. Then 
\begin{itemize}
\item[a)] The canonical quasiderivation $\Delta_{(S,\sigma)}$ associated to any choice of generating set $S$ and minimal presentation $\sigma$ is closable.
\item[b)] The domain ${\cal A}^1\Gamma$ of the canonical quasiderivation is independent of the choice of $S$ and $\sigma$. 
\item[c)] ${\cal A}^1\Gamma$ is an unconditional Banach algebra over $\C\Gamma$.
\item[d)] The canonical inclusion ${\cal A}^1\Gamma\hookrightarrow{\cal A}\Gamma$ is a homomorphism of Banach algebras and identifies ${\cal A}^1\Gamma$ with a dense and holomorphically closed subalgebra of ${\cal A}\Gamma$.
\end{itemize} 
\end{theorem}

\begin{proof}
\begin{itemize}
\item[a)] Let $(\alpha_n)$ be a sequence in ${\cal A}\Gamma$, such that
$$
\underset{n\to\infty}{\lim}(\alpha_n,\,\Delta_{(S,\sigma)}(\alpha_n))\,=\,
(0,\beta)\in {\cal A}\Gamma\times{\cal A}\Gamma\otimes_{uc}{\cal A}\Gamma.
$$ 
Consider the net $(\varphi_\alpha)$ of characteristic functions of finite 
subsets of $\Gamma$. The corresponding pointwise multiplication operators 
$M_{\varphi_\alpha}$ on ${\cal A}\Gamma$ are radial and satisfy 
$$\parallel M_{\varphi_\alpha}\parallel\leq 1,\,\forall\alpha,\;\text{and}\; 
\underset{\underset{\alpha}{\to}}{\lim}\,M_{\varphi_\alpha}x\,=\,x,\,\forall 
x\in {\cal A}\Gamma.$$ Consequently 
$M'_{\varphi_\alpha}=M_{\varphi_\alpha}\otimes_{uc}M_{\varphi_\alpha}$ 
defines a contracting operator on ${\cal A}\Gamma\otimes_{uc}{\cal A}\Gamma$
which satisfies 
$$\underset{\underset{n}{\to}}
{\lim}\,M'_{\varphi_\alpha}(\Delta_{(S,\sigma)}
(\alpha_n))\,=\,M'_{\varphi_\alpha}(\beta)\,=\,0,\forall\alpha,$$ and thus
$\beta\,=\,\underset{\underset{\alpha}{\to}}{\lim}\,M'_{\varphi_\alpha}
(\beta)\,=\,0$ as was to be shown.

\item[b)] Recall that a $(\lambda,C)$-quasigeodesic segment in a metric space $(X,d)$ is a map
$f:[a,b]\to X$ such that 
$$
\begin{array}{ccccc}
\frac{1}{\lambda}\vert s-t\vert-C & \leq & d(f(s),f(t)) & \leq & \lambda\vert s-t\vert+C \\
\end{array}
\eqno(3.9)
$$
for all $s,t\in[a,b]$. Suppose that $(X,d)$ is a $\delta$-hyperbolic space. Then a well known theorem (\cite{Gr}, 2.3, \cite{GH}, chapter 5,) asserts that there exists a constant $\mu(\lambda,C,\delta)$ such that the following holds: If $f:[a,b]\to X$ is a $(\lambda,C)$-quasigeodesic segment and $g:[a',b']\to X$ is a geodesic segment such that $f(a)=g(a')=x$ and $f(b)=g(b')=y$, then 
$$
f([a,b])\subset B(g([a',b']),\mu(\lambda,C,\delta)),
\eqno(3.10 )
$$
where $B(Y,r)$ denotes the $r$-neighbourhood of a subset $Y\subset X$ and the constant $\mu$ depends only on $\delta,\,\lambda,$ and $C$. Let now $\Gamma$ be a word-hyperbolic group and let $S,S'$ be finite symmetric sets of generators. Let $(X,d_S)$ be the Cayley graph of $\Gamma$ with respect to $S$, which is a $\delta$-hyperbolic geodesic metric space for suitable $\delta>0$. Let now $\sigma$ and $\sigma'$ be minimal presentations of the elements of $\Gamma$ with 
respect to $S$ and $S'$, respectively. Then there exist constants $\lambda', C'$ such that the following holds: for every $g\in\Gamma$ with $\sigma'(g)=(s_1',\ldots,s_n'),\,s_1',\ldots,s_n'\in S'$, the word 
$(\sigma(s_1'),\ldots,\sigma(s_n'))$ in the alphabet $S$ defines a $(\lambda',C')$-quasigeodesic in $(X,d_S)$. Thus the relation (3.10) implies by the reasoning used in the proof of 3.3 that
$$
\begin{array}{ccc}
\vert\Delta_{(S',\sigma')}(\alpha)\vert & \leq & \underset{l_S(t)\leq\mu(\lambda',C',\delta)}{\sum}\,i_t(\vert
\Delta_{(S,\sigma)}(\alpha)\vert)
\end{array}
\eqno(3.11)
$$
for every $\alpha\in\C\Gamma$. Thus 
$$
\begin{array}{ccc}
\parallel\Delta_{(S',\sigma')}(\alpha)\parallel_{{\cal A}\Gamma\otimes_{uc}{\cal A}\Gamma} & \leq & C''\parallel\Delta_{(S,\sigma)}(\alpha)\parallel_{{\cal A}\Gamma\otimes_{uc}{\cal A}\Gamma}
\end{array}
\eqno(3.12)
$$ 
for $C''\,=\,\underset{l_S(t)\leq\mu(\lambda',C',\delta)}{\sum}\parallel u_t\parallel\cdot\parallel u_{t^{-1}}\parallel$ and the assertion is proved.

\item[c)] Let $\Delta_{(S,\sigma)}$ be a canonical quasiderivation on $\C\Gamma$ and let $C_0>>0$ be a constant 
such that (3.3) holds. Then 
$$
\parallel\alpha\parallel'\,=\,C_0\parallel\alpha\parallel_{{\cal A}\Gamma}\,+\,
\parallel\Delta_{(S,\sigma)}(\alpha)\parallel_{{\cal A}\Gamma\otimes_{uc}{\cal A}\Gamma}
$$ 
is a submultiplicative unconditional seminorm on $\C\Gamma$ which is equivalent to the graph seminorm of $\Delta_{(S,\sigma)}$. The completion with respect to this seminorm is an unconditional Banach algebra over $\C\Gamma$ which coincides with ${\cal A}^1\Gamma$.

\item[d)] We have to show that the spectra of common elements of ${\cal A}^1\Gamma$ and ${\cal A}\Gamma$ coincide in both algebras. As ${\cal A}^1\Gamma$ is dense in ${\cal A}\Gamma$ it suffices to verify that 
$$
\begin{array}{cccc}
a\in{\cal A}^1\Gamma, & \parallel a\parallel_{{\cal A}\Gamma}\,<\,\frac{1}{C_0} & 
\Rightarrow & a\in GL_1({\cal A}^1\Gamma)
\end{array}
$$
where $C_0>>0$ is a constant satisfying (3.3). In fact, by (3.3)
$$
\parallel\Delta_{(S,\sigma)}(a^n)\parallel\,\leq\,
\underset{k=1}{\overset{n}{\sum}}\,C_0^k\parallel\Delta_{(S,\sigma)}(a)\parallel\cdot\parallel a\parallel_{{\cal A}\Gamma}^{n-1}
\,\leq\,C''\parallel\Delta_{(S,\sigma)}(a)\parallel\,(C_0\parallel a\parallel)_{{\cal A}\Gamma}^{n-1}
$$
which shows that the geometric series $\underset{k=0}{\overset{\infty}{\sum}}\,a^n$ converges in 
${\cal A}^1\Gamma$ against a\\ multiplicative inverse of $a$.
\end{itemize}
\end{proof}

\begin{cor}
Let $\Gamma$ be a word-hyperbolic group and let ${\cal A}\Gamma$ be an unconditional Banach algebra over the group ring $\C\Gamma$. Define inductively 
$$
\begin{array}{ccc}
{\cal A}^k\Gamma & = & ({\cal A}^{k-1})^1\Gamma.
\end{array}
\eqno(3.13)
$$
Then the intersection 
$$
\begin{array}{ccc}
{\cal A}^\infty\Gamma & = & \underset{k}{\bigcap}\,{\cal A}^k\Gamma
\end{array}
\eqno(3.14)
$$
is an admissible Fr\'echet algebra. It is unconditional over $\C\Gamma$ and dense and closed under holomorphic functional calculus in ${\cal A}\Gamma$. Any canonical quasiderivation $\Delta$ on $\C\Gamma$ extends to a bounded operator 
$$
\begin{array}{cccc}
\Delta: & {\cal A}^\infty\Gamma & \to & {\cal A}^\infty\Gamma\otimes_{uc}{\cal A}^\infty\Gamma
\end{array}
\eqno(3.15)
$$
\end{cor}
\begin{proof}
It follows from the previous theorem that ${\cal A}^k\Gamma$ is an unconditional Banach algebra over $\C\Gamma$ which is closed under holomorphic functional calculus in ${\cal A}\Gamma$. Thus ${\cal A}^\infty\Gamma$ is an unconditional Fr\'echet algebra over $\C\Gamma$ which is closed under holomorphic functional calculus in ${\cal A}\Gamma$ as well. We show that it is admissible, an open unit ball being given by 
$$
\begin{array}{ccc}
U & = & \{\alpha\in{\cal A}^\infty\Gamma,\,\parallel\alpha\parallel_{{\cal A}\Gamma}\,<\,1\}.
\end{array}
\eqno(3.16)
$$
So let $K\subset U$ be compact. As the multiplication in a Fr\'echet algebra is jointly continuous, each power $K^n,n\in\N,$ is compact and it remains to show that their union is still relatively compact. For this it suffices to verify that $\underset{n\to\infty}{\lim}\parallel K^n\parallel_{{\cal A}^k\Gamma}=0$ for all $k\in\N$. We argue by induction on $k$. Choose a generating set $S\subset\Gamma$ and a minimal presentation $\sigma$ of the elements of $\Gamma$ and let $\Delta_{(S,\sigma)}$ be the associated quasiderivation. Any Banach norm on ${\cal A}^1\Gamma$ is then equivalent to the norm 
$\parallel\alpha\parallel'=\parallel\alpha\parallel_{{\cal A}\Gamma}+\parallel\Delta_{(S,\sigma)}(\alpha)\parallel_{uc}.$ Let $C_0>>0$ be a constant such that (3.3) holds. Then \\ $\parallel\Delta_{(S,\sigma)}(K^{2^n})\parallel_{uc}\,\leq\,(2C_0)^n\,\parallel\Delta_{(S,\sigma)}(K)\parallel_{uc}\,\cdot\parallel K\parallel_{{\cal A}\Gamma}^{2^n-1}$ by a simple induction and more generally 
$\parallel\Delta_{(S,\sigma)}(K^n)\parallel_{uc}\,\leq\,(2C_0)^{log_2(n)}\,\parallel\Delta_{(S,\sigma)}(K)\parallel_{uc}\,\cdot\parallel K\parallel_{{\cal A}\Gamma}^{n-1}$. As $\parallel K\parallel_{{\cal A}\Gamma}<1$ by our choice of $U$ this implies $\underset{n\to\infty}{\lim}\parallel K^n\parallel_{{\cal A}\Gamma}=0$ and\\ $\underset{n\to\infty}{\lim}\parallel\Delta_{(S,\sigma)}(K^n)\parallel_{uc}=0$ which proves the first induction step. Iterating the argument yields our claim for every $k\in\N$.

Finally, if $\Delta_{(S',\sigma')}$ is any canonical quasiderivation on $\C\Gamma$, then, according to (3.12),\\ $\parallel\Delta_{(S',\sigma')}(\alpha)\parallel_{{\cal A}\Gamma\otimes_{uc}{\cal A}\Gamma}\leq  C''\parallel\Delta_{(S,\sigma)}(\alpha)\parallel_{{\cal A}\Gamma\otimes_{uc}{\cal A}\Gamma}$ for some constant $C''>>0$. This shows $\parallel\Delta_{(S',\sigma')}(\alpha)
\parallel_{{\cal A}^k\Gamma\otimes_{uc}{\cal A}^k\Gamma}\,\leq\,
C''\parallel\alpha\parallel_{k+1}$ and finishes the proof of the corollary.
\end{proof} 
\begin{remark}
Conclusions analogous to 3.5 and 3.6 hold also under the weaker assumption that ${\cal A}\Gamma$ is an unconditional admissible Fr\'echet algebra over $\C\Gamma$. The demonstrations, which are similar
to the case of Banach algebras treated above, have been omitted as notations become much more elaborate in the general case.
\end{remark}

\section{Minimal subalgebras and Sobolev norms}

It turns out that the subalgebra ${\cal A}^\infty\Gamma$ of infinitely quasiderivable elements of 
${\cal A}\Gamma$ is not yet small enough for our purpose. Every quasidifferential operator of finite order is bounded on ${\cal A}^\infty\Gamma$ but we will need to control the norm of quasidifferential operators of infinite order. To this end we introduce subalgebras of ${\cal A}\Gamma$ which correspond to algebras of analytic functions.

\begin{definition}
Let $\Gamma$ be a discrete group and let $S$ be a finite symmetric set of generators.
Let ${\cal A}\Gamma$ be an unconditional admissible Fr\'echet algebra over $\C\Gamma$ and let $U\subset{\cal A}\Gamma$ be an open unit ball which is stable under $a\mapsto\vert a\vert$. Put

$$
T_n\,=\,(\underset{l_S(g)\leq n}{\sum}\C\,u_g)\,\bigcap\,\frac{n}{n+1}\cdot U
\eqno(4.1)
$$ 

We denote by $\parallel-\parallel_n$ the largest unconditional submultiplicative seminorm on $\C\Gamma$ (see 2.11) such that
$$
\begin{array}{ccc}
\parallel T_n\parallel_n & \leq & 1.
\end{array}
\eqno(4.2)
$$
The completion of $\C\Gamma$ with respect to this norm is an unconditional Banach algebra over $\C\Gamma$ denoted by ${\cal A}_n\Gamma$. For $n<m$ the identity on $\C\Gamma$ extends to a bounded homomorphism ${\cal A}_n\Gamma\to{\cal A}_{m}\Gamma$. The inductive system 
$$
\begin{array}{ccc}
"{\cal A}_\omega\Gamma" & = & "\underset{\underset{n}{\rightarrow}}{\lim}"\,{\cal A}_n\Gamma \\
\end{array}
\eqno(4.3)
$$
is independent of the choice of $S\subset\Gamma$ and $U\subset{\cal A}\Gamma$.
The union
$$
\begin{array}{ccc}
{\cal A}_\omega\Gamma & = & \underset{n}{\bigcup}\,{\cal A}_n\Gamma \\
\end{array}
\eqno(4.4)
$$
is a bornological algebra which is independent of the choice of $S\subset\Gamma$ and $U\subset{\cal A}\Gamma$.
\end{definition}

\begin{lemma}
In the notations of 4.1 one has 
$$
Sp_{{\cal A}_{n+1}\Gamma}(\alpha)\subset Sp_{{\cal A}_n\Gamma}(\alpha)
\eqno(4.5)
$$ 
and 
$$
\begin{array}{ccc}
Sp_{{\cal A}\Gamma}(\alpha) & = & \underset{n}{\bigcap}\,Sp_{{\cal A}_n\Gamma}(\alpha)
\end{array}
\eqno(4.6)
$$
for every $\alpha\in{\cal A}_\omega\Gamma$.
\end{lemma}
This is clear because ${\cal A}\Gamma$ is admissible and contains $\C\Gamma$ as dense subalgebra.
We define now Sobolev norms on the minimal Banach subalgebras introduced above.

\begin{definition}
Let the notations of 4.1 and 4.2 be understood. 
\begin{itemize}
\item[a)] Denote by $\parallel-\parallel_{n,k},\,n,k\in\N,$ the largest 
unconditional seminorm on $\C\Gamma$ such that
$$
\begin{array}{cccc}
\parallel T_n^m\parallel_{n,k} & \leq & m^k, & \forall m\in\N.
\end{array}
\eqno(4.7)
$$
\item[b)] Denote by the same symbol the largest unconditional seminorm on $\C\Gamma\otimes\C\Gamma$ such that 
$$
\begin{array}{cccc}
\parallel T_n^{m_1}\otimes T_n^{m_2}\parallel_{n,k} & \leq & (m_1+m_2)^k, & \forall m_1,\,m_2\in\N.
\end{array}
\eqno(4.8)
$$
\end{itemize}
\end{definition}
In particular, the multiplication map $m:\C\Gamma\otimes\C\Gamma\to\C\Gamma$ satisfies 
$$
\begin{array}{ccc}
\parallel m(\alpha)\parallel_{n,k} & \leq & \parallel\alpha\parallel_{n,k}
\end{array}
\eqno(4.9)
$$ 
for all $\alpha\in\C\Gamma\otimes\C\Gamma$ and all $n,k\in\N$.
We estimate the Sobolev norms of the canonical quasiderivation.
\begin{lemma}
Let $\Gamma$ be a word-hyperbolic group. Let $S$ be a finite symmetric set of generators and let $\sigma$ be a minimal presentation of the elements of $\Gamma$ in terms of the alphabet $S$. 
Suppose that $\delta\geq 0$ is such that the Cayley graph ${\cal G}(\Gamma,S)$ is a $\delta$-hyperbolic geodesic metric space. 

Let ${\cal A}\Gamma$ be an unconditional admissible Fr\'echet algebra over $\C\Gamma$ and choose an open unit ball $U\subset{\cal A}\Gamma$ which is stable under $a\mapsto\vert a\vert$. 

Consider the associated minimal Banachalgebras ${\cal A}_n\Gamma,\,n\in\N,$ (see 4.1) and let \\ $\parallel-\parallel_{n,k},\,n,k\in\N,$ be the corresponding Sobolevnorms of 4.3. 

Then there exist constants $C_1(n)=C_1(n,S,\delta)>0,\,C_2=C_2(S,\delta)>0,$ and \\ $d=d(S,\delta)\in\N_+,$ such that
the canonical quasiderivation $\Delta_{(S,\sigma)}$ (3.1) satisfies
$$
\begin{array}{ccc}
\parallel\Delta_{(S,\sigma)}(\alpha)\parallel_{n,k} & \leq & C_1(n)
\cdot C_2^k\cdot\parallel\alpha\parallel_{n,k+d}
\end{array}
\eqno(4.10)
$$
for all $\alpha\in\C\Gamma,$ and all $n,k\in\N$.
\end{lemma}

\begin{proof}
Let $\alpha_1,\ldots,\alpha_m\in T_n$ (note that $T_n$ is stable under $a\mapsto\vert a\vert$.) 
From (3.5) we derive by induction 
$$
\begin{array}{ccc}
\vert\Delta_{(S,\sigma)}(\alpha_1\ldots\alpha_m)\vert & \leq & 
\underset{j=1}{\overset{l}{\sum}}\,\underset{t_1,...,t_j}{\sum}\,\underset{i=1}{\overset{m}{\sum}}
\vert\alpha_1\vert\ldots\vert\alpha_{i-1}\vert\,(i_{t_1}\ldots i_{t_j}\vert\Delta_{(S,\sigma)}(\alpha_i)\vert)
\vert\alpha_{i+1}\vert\ldots\vert\alpha_m\vert
\end{array}
$$
 where the sum is taken over all $t_1,...,t_j$ such that $l_S(t_i)\leq\delta,\,i=1,\ldots,j,$ and $l=[\log_2(m)]+1$.

As $T_n$ is contained in the linear span of the finite set $\{u_g,\,l_S(g)\leq n\}$, there exists a constant $C_3(n)>>0$ and  finite families $\alpha_q',\alpha_q''\in T_n,\,q\in I,$ such that the following holds:
for very $\alpha\in T_n$ there are numbers $\lambda_q\geq 0,\,q\in I,$ such that \\ $\vert\Delta_{(S,\sigma)}(\alpha)\vert\leq\underset{I}{\sum}
\lambda_q\vert\alpha_q'\vert\otimes\vert\alpha_q''\vert$ and $\underset{I}{\sum}\lambda_q\leq C_3(n)$.

From this we deduce
$$
\vert\alpha_1\vert\ldots\vert\alpha_{i-1}\vert\,(i_{t_1}\ldots i_{t_j}\vert\Delta_{(S,\sigma)}(\alpha_i)\vert)
\vert\alpha_{i+1}\vert\ldots\vert\alpha_m\vert
$$
$$
\leq\underset{I}{\sum}\lambda_q
\vert\alpha_1\vert\ldots\vert\alpha_{i-1}\vert\,(i_{t_1}\ldots i_{t_j}\vert\alpha_q'\vert\otimes\vert\alpha_q''\vert)
\vert\alpha_{i+1}\vert\ldots\vert\alpha_m\vert
$$
$$
=\underset{I}{\sum}\lambda_q
\vert\alpha_1\vert\ldots\vert\alpha_{i-1}\vert\vert\alpha_q'\vert u_{t_j}^{-1}\ldots u_{t_1}^{-1}\otimes
u_{t_1}\ldots u_{t_j}\vert\alpha_q''\vert
\vert\alpha_{i+1}\vert\ldots\vert\alpha_m\vert
$$
and thus 
$$
\parallel\vert\alpha_1\vert\ldots\vert\alpha_{i-1}\vert\,(i_{t_1}\ldots i_{t_j}\vert\Delta_{(S,\sigma)}(\alpha_i)\vert)
\vert\alpha_{i+1}\vert\ldots\vert\alpha_m\vert\parallel_{n,k}
$$
$$
\leq\underset{I}{\sum}\lambda_q\left(i+j\delta+j\delta+m-i+1\right)^k
C_4^{2j}\,\leq\,
C_3(n)(m+1+2j\delta)^k\cdot C_4^{2j}
$$
where $C_4\,=\,\Inf\{C^{-1},\,2Cu_t\in U,\,t\in\Gamma,\,l_S(t)\leq\delta\}$.\\
Choose a constant $C_5(\delta)>0$ such that 
$$
m+1+2([\log_2(m)]+1)\delta\leq C_5(\delta)m.
$$ 
for all $m\geq 1$. As $j\leq [\log_2(m)]+1$ the previous estimates imply 
$$
\parallel\vert\alpha_1\vert\ldots\vert\alpha_{i-1}\vert\,(i_{t_1}\ldots i_{t_j}\vert\Delta_{(S,\sigma)}(\alpha_i)\vert)
\vert\alpha_{i+1}\vert\ldots\vert\alpha_m\vert\parallel_{n,k}\,\leq\,
C(n,\delta)\cdot C_5(\delta)^k\cdot m^{k+C_6}
$$
Summing up, one finds
$$
\parallel\Delta_{(S,\sigma)}(\alpha_1\ldots\alpha_m)\parallel_{n,k}\, \leq\,
C(n,\delta)\cdot C_5(\delta)^k\underset{j=1}{\overset{[\log_2(m)]+1}{\sum}}\,
(\sharp\{ t\in\Gamma,\,l_S(t)\leq\delta\})^j\cdot m^{k+C_6+1}
$$
$$
\leq\,C'(n,\delta)\cdot C_5(\delta)^k\cdot m^{k+C_7(\delta)}
$$
for a suitable integer $C_7(\delta)$. As the map $\Delta_{(S,\sigma)}$ is radial, this shows in view of the definition of the Sobolev norms that 
$$
\parallel\Delta_{(S,\sigma)}(\alpha)\parallel_{n,k}\,\leq\,C'(n,\delta)\cdot C_5(\delta)^k\cdot
\parallel\alpha\parallel_{n,k+C_7(\delta)}
$$
for all $\alpha\in\C\Gamma$
which proves our claim.
\end{proof}

\begin{definition}
Let $\Gamma$ be a group with finite symmetric set $S$ of generators and
let $\sigma$ be a minimal presentation of the elements of $\Gamma$ in terms of the alphabet $S$. 
 A radial map $\psi:\C\Gamma\to\C\Gamma\otimes\C\Gamma$ is called $R$-special, $R>0$, if
\begin{itemize}
\item[i)] 
$$
Supp(\psi(u_g))\subset m^{-1}(g), \forall g\in\Gamma
$$ 
where $m:\Gamma\times\Gamma\to\Gamma$ is the group multiplication.
\item[ii)] 
$$
\underset{g\in\Gamma}{\Sup}\parallel\psi(u_g)\parallel_{\ell^\infty(\Gamma\times\Gamma)}\leq 1.
$$
\item[iii)] 
$$
\pi_1(Supp(\psi(u_g)))\subset B(\sigma(g),R)\subset\Gamma,\, \forall g\in\Gamma
\,=\,{\cal G}(\Gamma,S)_0,
$$
where $B(-,R)$ denotes the $R$-neighbourhood of a given subset in the Cayley graph of $(\Gamma,S)$.
\end{itemize}
\end{definition}

\begin{prop}
Let $\Gamma$ be a word-hyperbolic group. Let $S$ be a finite symmetric set of generators and let $\sigma$ be a minimal presentation of the elements of $\Gamma$ in terms of the alphabet $S$. 
Suppose that $\delta\geq 0$ is such that the Cayley graph ${\cal G}(\Gamma,S)$ is a $\delta$-hyperbolic geodesic metric space. 

Let ${\cal A}\Gamma$ be an unconditional admissible Fr\'echet algebra over $\C\Gamma$ and choose an open unit ball $U\subset{\cal A}\Gamma$ which is stable under $a\mapsto\vert a\vert$. 

Consider the associated minimal Banachalgebras ${\cal A}_n\Gamma,\,n\in\N,$
and let $\psi:\C\Gamma\to\C\Gamma\otimes\C\Gamma$ be an $R$-special map.
Then there exist constants\\ $C_8(n)=C_8(n,R,S,\delta)>0,\,C_9=C_9(R,S,\delta)>0,$ and $d=d(S,\delta)\in\N_+,$ such that
$$
\begin{array}{ccc}
\parallel\psi(\alpha)\parallel_{n,k} & \leq & C_8(n,R)
\cdot C_9(R)^k\cdot\parallel\alpha\parallel_{n,k+d}
\end{array}
\eqno(4.11)
$$
for all $\alpha\in\C\Gamma,$ where $\parallel-\parallel_{n,k},\,n,k\in\N,$ denotes the Sobolev norm introduced in 4.3.
\end{prop}

\begin{proof}
By definition of a special map one has
$$
\begin{array}{ccc}
\vert\psi(\alpha)\vert & \leq &
\underset{l_S(t)\leq R}{\sum}\,i_t(\vert\Delta_{(S,\sigma)}(\alpha)\vert)
\end{array}
\eqno(4.12)
$$
for all $\alpha\in\C\Gamma$. The assertion follows then from 4.4 and the estimates used in its proof.
\end{proof}

\section{Estimates on conjugacy classes}

We recall Gromov's solution of the {\bf conjugacy problem} in hyperbolic groups.

\begin{theorem} (Gromov), see \cite{Gr}, 7.4.\\
Let $\Gamma$ be a word-hyperbolic group with finite symmetric set of generators $S$. There exists a constant $C_{10}=C_{10}(\Gamma,S)>>0$ such that the following holds: if $g,g'\in\Gamma$ are conjugate elements, then there exists $u\in\Gamma$ such that $g'=ugu^{-1}$ and 
$$
\begin{array}{ccc}
l_S(u) & \leq & C_{10}\left(l_S(g)+l_S(g')\right) \\
\end{array}
\eqno(5.1)
$$
\end{theorem}

Let $\Gamma$ be a word-hyperbolic group with finite symmetric set of generators $S$. Let $\sigma$ be a minimal presentation of the elements of $\Gamma$ by words in the alphabet $S$. For an element $x\in\Gamma$ denote by $\langle x\rangle\subset\Gamma$ its conjugacy class. We choose in every conjugacy class an element $h=h(\langle x\rangle)$ of minimal word length. Furthermore, we choose for every nontrivial element $g\in\Gamma$ an element 
$u=u(g)$ of minimal word length in
$$
S(g)\,=\,\{v\in\Gamma, v^{-1}gv=h(\langle g\rangle)'\}.
$$
where $h(\langle g\rangle)'$ is an element of $\Gamma$ which can be represented by a cyclic permutation of the word $\sigma(h(\langle g\rangle))$ for $h(\langle g\rangle)$.
Note that 
$$
l_S(u(g))\,\leq\,2\,C_{10}\,l_S(g)
\eqno(5.2)
$$ 
by Gromov's theorem. We suppose henceforth that $C_{10}\geq 1$.

It is obvious from the definitions that every point of an edge of a geodesic rectangle in a $\delta$-hyperbolic geodesic metric space, $\delta\in\N_+$, is at distance at most $2\delta$ from the union of the three other edges. 

For $g\in\Gamma$ we consider the geodesic rectangle in the Cayley graph $\cal{G}(\Gamma,S)$ with corners 
$$
A=e,\,B=u(g),\,C=u(g)h(\langle g\rangle)',\,D=g,
$$ 
and edges 
$$
a=\sigma(u(g)),\,b=\,\text{cyclic permutation of}\,\sigma(h(\langle g\rangle)), c=\sigma(u(g)^{-1}),\,d=\sigma(g).
$$
Define now a map $\mu:\Gamma\to\Gamma$ as follows:
\begin{itemize}
\item[i)] $\mu(e)=e$ \\
\item[ii)] Otherwise if $u(g)=e$ and $\sigma(h(\langle g\rangle))=s_1\ldots s_n$, then 
$\mu(g)=s_{m+1}\ldots s_n$ where $0\leq m\leq n$ is the largest integer such that 
$s_{m+1}\ldots s_ns_1\ldots s_m=g$.
\item[iii)] Otherwise $\mu(g)=u(g)$ if $l_S(u(g))\leq 2\delta$ \\
\item[iv)] Otherwise $\mu(g)$ is the product of the first $2\delta$ letters of $\sigma(u(g))$ if\\ 
$l_S(u(g))\leq 24 C_{10}\delta$ where $C_{10}$ is Gromov's constant (5.1).\\
\item[v)] Otherwise $\mu(g)$ is some vertex of ${\cal G}(\Gamma,S)$ situated on the edge $b$ at distance 
at most $2\delta+\frac12$ from the midpoint of $d$ if such a vertex exists.\\
\item[vi)] Otherwise $\mu(g)$ is some vertex of ${\cal G}(\Gamma,S)$ situated on one of the edges $a$ or $c$ at distance at most $2\delta+\frac12$ from the midpoint of the edge $d$.\\
\end{itemize}
Put
$$
\begin{array}{cc}
\varphi:\Gamma\to\Gamma, & g\mapsto \mu(g)^{-1}g\mu(g),
\end{array}
\eqno(5.3)
$$ 
We are interested in the limit
$$
\begin{array}{ccc}
\Phi & = & \underset{n\to\infty}{\lim}\,\varphi^n
\end{array}
\eqno(5.4)
$$
and analyse in particular how fast it is attained.
\begin{lemma}
There exist constants $C_{11},C_{12}>>0$ 
such that 
$$
\varphi^n(g)\,=\,h(\langle g\rangle)
\eqno(5.5)
$$ 
provided that $g\in\Gamma$ satisfies 
$$
C_{11}\,log(l_S(g))\,+\,C_{12}\,\leq\,n.
\eqno(5.6)
$$
\end{lemma}

\begin{proof}
The definition of the map $\mu$ (and thus of $\varphi$) distinguishes six separate cases which have to be treated one after the other. 

Case i): there is nothing to show. \\

Case ii): $\varphi(g)\,=\,\varphi^j(g)\,=\,
h(\langle g\rangle)$ for all $j\geq 1$. \\

Case iii): $\varphi(g)$ satisfies the condition ii) so that $\varphi^2(g)=h(\langle g\rangle)$. \\

Case iv):
$$
l_S(u(\varphi(g)))\,\leq\,l_S(u(g))-2\delta\,\leq\,24\,C_{10}\,\delta
$$ so that $\varphi^j(g)$ falls 
into case ii) or iii) provided that $j\,\geq\,12\,C_{10}$. So 
$$
\varphi^m(g)=h(\langle g\rangle)\,\,\,\text{if}\,\,\,m\,\geq\,12\,C_{10}\,+\,2.
$$ 
Case v): $\varphi(g)$ may be represented by a cyclic permutation of the word for $h(\langle g\rangle)$ so that $\varphi^2(g)=h(\langle g\rangle)$ according to case ii). \\

Case vi): Let $g\in\Gamma$. We calculate an upper bound for the smallest integer $k$ such that one of the conditions i) to v) applies to $\varphi^k(g)$. Suppose that $\mu(g)$ is a vertex of the Cayley graph situated on the edge $a=\sigma(u(g))$ with endpoints $e$ and $u(g)$ at distance at most $2\delta+\frac12<3\delta$ from the midpoint $m'$ of $d=[e,g]$ (the calculation in the case $\mu(g)\in c$ is identical). Note that, according to Gromov's theorem (5.2) $l_S(u(g))\leq 2\,C_{10}\,l_S(g)$ or $l_S(g)\geq\frac{1}{2C_{10}}l_S(u(g))$. Moreover 
$l_S(u(g))> 24 C_{10}\delta$ because we are not in case iv). One finds 
$$
l_S(\mu(g))\,=\,d_S(e,\mu(g))\,\geq\,
d_S(e,m')-d_S(\mu(g),m')\,\geq\,\frac12\,l_S(g)-3\delta
$$
$$
\geq\,\frac{1}{4C_{10}}l_S(u(g))-
\frac{1}{8C_{10}}l_S(u(g))\,=\,\frac{1}{8C_{10}}l_S(u(g))
$$
Therefore 
$$
\varphi(g)\,=\,\mu(g)^{-1}g\mu(g)\,=\,(\mu(g)^{-1}u(g))h(\langle g\rangle)'(u(g)^{-1}\mu(g))
$$ 
and  
$$
l_S(u(\varphi(g)))\,\leq\,l_S(\mu(g)^{-1}u(g))\,=\,d_S(\mu(g),u(g))
$$
$$
=\,l_S(u(g))-l_S(\mu(g))\,\leq\,C_{13}\,l_S(u(g))
$$
where
$$
0\,<\,C_{13}\,=\,1-\frac{1}{8C_{10}}\,<\,1.
$$
Suppose now that condition vi) above applies to the elements $g,\varphi(g),\ldots,\varphi^{j}(g)$. Then, by 
our calculation, 
$$
24\,C_{10}\,\delta\,<\,l_S(u(\varphi^j(g)))\,\leq\,C_{13}^j\,l_S(u(g))\,\leq\,2\,C_{10}\,C_{13}^j\,l_S(g)
$$
which yields
$$
j\,\leq\,\log(C_{13}^{-1})^{-1}\,\log(l_S(g)).
$$
So for a suitable constant $C_{14}$ the element $\varphi^m(g),\,m\,=\,[C_{14}\log(l_S(g))]+1$ falls in one of the first five cases. This completes the proof of the lemma. 
\end{proof}

\begin{lemma}
The notations made so far in the present section are understood. The radial linear map 
$$
\begin{array}{ccc}
\varphi: & \C\Gamma\to\C\Gamma, & u_g\mapsto u_{\varphi(g)}
\end{array}
\eqno(5.7)
$$ 
can be written as the composition 
$$
\begin{array}{ccc}
\varphi & = & m\,\circ\,sw\,\circ\,\psi
\end{array}
\eqno(5.8)
$$
where
$$
\begin{array}{ccc}
\psi: & \C\Gamma\to\C\Gamma\otimes\C\Gamma, & u_g\mapsto u_{\mu(g)}\otimes u_{\mu(g)^{-1}g}
\end{array}
\eqno(5.9)
$$
is radial and $2\delta+1$-special in the sense of 3.5,
$sw:\,\C\Gamma\otimes\C\Gamma\to\C\Gamma\otimes\C\Gamma$
is the radial linear map which flips the factors, and 
$m:\C\Gamma\otimes\C\Gamma\to\C\Gamma$
is the multiplication, which is radial as well.
\end{lemma}

 We are now ready to state our main result.

\begin{theorem}
Let $\Gamma$ be a word-hyperbolic group. Let ${\cal A}\Gamma$ be an unconditional admissible Fr\'echet algebra over $\C\Gamma$ and let ${\cal A}_\omega\Gamma$ be its canonical, holomorphically closed bornological subalgebra (4.4). Let $S$ be a finite symmetric set of generators of $\Gamma$ and choose an element $h(\langle x\rangle)$ of minimal word length (with respect to $S$) in every conjugacy class $\langle x\rangle$ of $\Gamma$. Then the characteristic linear map 
$$
\begin{array}{cc}
 \C\Gamma\to\C\Gamma, & u_g\mapsto u_{h(\langle g\rangle)}
\end{array}
\eqno(5.10)
$$
extends to a bounded linear operator
$$
\begin{array}{cc}
{\bf \Phi}: & {\cal A}_\omega\Gamma\,\to\,{\cal A}_\omega\Gamma
\end{array}
\eqno(5.11)
$$
To be more precise, if ${\cal A}_\omega\Gamma\,=\,\underset{n}{\bigcup}\,{\cal A}_n\Gamma$ is viewed as monotone union of the minimal Banach subalgebras ${\cal A}_n\Gamma\subset{\cal A}\Gamma$, then for $n<m$ the characteristic linear map extends to a bounded linear operator
$$
\begin{array}{cc}
{\bf \Phi}: & {\cal A}_n\Gamma\,\to\,{\cal A}_m\Gamma.
\end{array}
\eqno(5.12)
$$
\end{theorem}

\begin{proof}
Clearly (5.12) implies (5.11). (Note that, contrary to the 
union ${\cal A}_\omega\Gamma$, the individual Banach algebras ${\cal A}_n\Gamma$ depend on the choice of a generating set $S$ of $\Gamma$ and of an open unit ball $U$ of ${\cal A}\Gamma$, which we may assume to be stable under $a\mapsto\vert a\vert$.) Recall the definition of ${\cal A}_n\Gamma$. 
If $T_n\,=\,(\underset{l_S(g)\leq n}{\sum}\C\,u_g)\,\bigcap\,\frac{n}{n+1}\cdot U$ then ${\cal A}_n\Gamma$ is the completion of $\C\Gamma$ with respect to the largest unconditional submultiplicative seminorm $\parallel-\parallel_n$ such that $\parallel T_n\parallel_n\leq 1$. As $T_n\subset T_{n+1}$ the structure homomorphisms ${\cal A}_n\Gamma\to{\cal A}_{n+1}\Gamma$ are contractive.
It suffices to show that for every $n\in\N$ there exists a constant $C=C(n)>>0$ such that 
$$
\parallel\Phi(\alpha)\parallel_{n+1}\,\leq\,C\parallel\alpha\parallel_n
$$ 
for all $\alpha\in\C\Gamma$.
The characteristic linear map $\Phi$ is radial by construction so that the pullback seminorm $\Phi^*(\parallel-\parallel_{n+1})$ is still an unconditional norm on $\C\Gamma$. On the other hand $\parallel-\parallel_n$ is the largest unconditional seminorm $\parallel-\parallel'$ such that $\parallel T_n^N\parallel'\leq 1$ for all $N$ (see 2.11).
Our claim follows thus for given $n$ from the estimate
$$
\begin{array}{ccc}
\underset{N}{\Sup}\,\parallel\Phi(T_n^N)\parallel_{n+1} & < & \infty
\end{array}
\eqno(5.13)
$$
which we will establish now.
The set $T_n\subset\C\Gamma$ is contained in the linear span of group elements of word length at most $n$ so that $T_n^N$ consists of linear combinations of group elements of word length at most $Nn$. From lemma 5.2 we deduce therefore 
$$
\Phi(T_n^N)\,=\,\varphi^m(T_n^N)
\eqno(5.14)
$$
with 
$$
m\,\leq\,C_{15}(n)\cdot\log(N)\,+\,C_{16}(n).
\eqno(5.15)
$$
Recall (see 5.3) that the operator $\varphi$ may be viewed as composition $\varphi=m\circ sw\circ\psi$ of radial linear maps. This leads to the following estimates for the Sobolev-norms 4.3, which make use of 4.3, 4.6, and of the fact that $\psi$ is $2\delta+1$ special:
$$
\parallel\varphi(\alpha)\parallel_{n+1,k}\,\leq\,C_8(n+1)\,C_9^k\parallel\alpha\parallel_{n+1,k+d}
\eqno(5.16)
$$
for suitable $C_8(n+1)=C_8(n+1,S,\delta),\,C_9(S,\delta)$ and $d=d(S,\delta)$. 
Iteration yields
$$
\parallel\Phi(T_n^N)\parallel_{n+1}\,=\,\parallel\varphi^m(T_n^N)\parallel_{n+1,0}
\,\leq\,C_8^m\,C_9^{(1+\ldots+(m-1))d}\parallel T_n^N\parallel_{n+1,md}
\eqno(5.17)
$$
For the latter norm one derives from
$$
\begin{array}{ccc}
T_n & \subset & \frac{n}{n+1}\cdot\frac{n+2}{n+1}\cdot T_{n+1} 
\end{array}
\eqno(5.18)
$$
the estimate
$$
\parallel T_n^N\parallel_{n+1,md}\,\leq\,C_{17}(n)^N\parallel T_{n+1}^N\parallel_{n+1,md}\,
\leq\,C_{17}(n)^N\,N^{md}
\eqno(5.19)
$$
where 
$$
0<C_{17}(n)=\frac{n(n+2)}{(n+1)^2}<1
\eqno(5.20)
$$
and finds
$$
\parallel\Phi(T_n^N)\parallel_{n+1}\,\leq\,C_8^m\,C_9^{(1+\ldots+(m-1))d}C_{17}(n)^N\,N^{md}\,
\leq\,C_{17}^N\,C_{18}^{m^2}\,N^{md}.
\eqno(5.21)
$$
Taking logarithms we obtain finally
$$
\log(\parallel\Phi(T_n^N)\parallel_{n+1})\,\leq\,m^2\,\log(C_{18})\,+\,md\log(N)-N\log(C_{17}^{-1})
$$
$$
\leq\,C_{19}(\log(N))^2\,+\,C_{20}\log(N)\,-\,N\log(C_{17}^{-1})
\eqno(5.22)
$$
As $\log(C_{17}^{-1})>0$ by (5.20) the last term is bounded from above independently of $N$. 
This establishes (5.13) and therefore the theorem.
\end{proof}

\section{Consequences}

\subsection{Traces on unconditional completions of group rings}

Recall that a class function on a group $\Gamma$ is a function which is constant on conjugacy classes. The class functions correspond bijectively to the traces on the group ring under the linear isomorphism between the space of complex valued functions on $\Gamma$ and the dual space of the group ring $\C\Gamma$.

\begin{definition}
A class function $\tau$ on a finitely generated group $\Gamma$ is called \\ {\bf tempered} if 
$$
\begin{array}{ccccc}
\underset{\langle x\rangle}{\sum}\,(1+l_S(\langle x\rangle))^{-2k}\,\vert\tau(\langle x\rangle)\vert^2
& < & \infty & \text{for} & k>>0,
\end{array}
\eqno(6.1)
$$
where $l_S(\langle x\rangle)$ denotes the minimal word length (with respect to a finite symmetric generating set $S$) of an element in the conjugacy class $\langle x\rangle$. 
\end{definition}

Note that this condition does not depend on the choice of $S$.
 In particular, every class function which is supported in only finitely many conjugacy classes is tempered.

\begin{theorem}

Let $\Gamma$ be a word-hyperbolic group and let $C^*_r(\Gamma)$ be its reduced group $C^*$-algebra.
There exists an infinite increasing sequence 
$({\cal A}_n\Gamma)$ of Banach subalgebras of $C^*_r(\Gamma)$, containing $\C\Gamma$, and such that ${\cal A}_\omega\Gamma\,=\,\underset{n}{\bigcup}\,{\cal A}_n\Gamma$ has the following properties:
\begin{itemize}
\item[$\bullet$] ${\cal A}_\omega\Gamma$
is dense and closed under holomorphic functional calculus in $C^*_r(\Gamma)$. \\
\item[$\bullet$] Every tempered class function on $\Gamma$ gives rise to a bounded trace on 
${\cal A}_\omega\Gamma$.
\end{itemize}
\end{theorem}

\begin{proof}
Choose a finite symmetric set of generators $S$ and let $l_S$ the corresponding word length function.
Denote by ${\mathfrak A}\Gamma$ the space of square-summable functions of rapid decay on $\Gamma$, i.e. the completion of the group ring $\C\Gamma$ with respect to the unconditional seminorms
$$
\begin{array}{cccc}
\parallel a\parallel_k^2 & = & \underset{g}{\sum}\,(1+l_S(g))^{2k}\,\vert a_g\vert^2, & k\in\N. \\
\end{array}
\eqno(6.2)
$$
It was shown by Haagerup for free groups \cite{Ha} and lateron by Jolissaint \cite{Jol} for general word-hyperbolic groups, that ${\mathfrak A}\Gamma$ is actually a Fr\'echetalgebra under convolution and is closed under holomorphic functional calculus in $C^*_r(\Gamma)$. Moreover it
 is admissible \cite{Pu3}, 4.2. 

Let ${\mathfrak A}_n\Gamma,\,n\in\N,$ be the minimal Banach subalgebras 4.1 of ${\mathfrak A}\Gamma$ and let ${\mathfrak A}_\omega\Gamma$ be their union. Then ${\mathfrak A}_\omega\Gamma$ is closed under holomorphic functional calculus in ${\mathfrak A}\Gamma$ and thus in $C^*_r(\Gamma)$.

Choose in every conjugacy class $\langle x\rangle$ of $\Gamma$ an element $h(\langle x\rangle)$ of minimal word length and denote by $\Phi$ the associated characteristic map (5.10), (5.11). Let $\tau$ be a class function on $\Gamma$. Then 
$a\mapsto\underset{\langle x\rangle}{\sum}\tau(\langle x\rangle)a_{h(\langle x\rangle)}$ defines a linear functional $\ell_\tau$ on $\C\Gamma$ such that the composition $\tau'=\ell_\tau\circ\Phi$ is the trace on $\C\Gamma$ associated to $\tau$.

According to theorem 5.4, the characteristic map $\Phi$ extends to a bounded endomorphism of  ${\mathfrak A}_\omega\Gamma$. A look at the definitions shows that $\ell_\tau$ extends to a bounded linear functional on the Jolissaint algebra ${\mathfrak A}\Gamma$ iff $\tau$ is tempered in the sense of 6.1. Consequently the trace $\tau'=\ell_\tau\circ\Phi$ extends to a bounded linear functional on ${\mathfrak A}_\omega\Gamma$ in this case as well.
\end{proof}

\begin{cor}
Let $\Gamma$ be a word-hyperbolic group and let ${\cal A}\Gamma$ be an unconditional admissible Fr\'echet-algebra over $\C\Gamma$. Denote by ${\cal A}_\omega\Gamma$ the canonical holomorphically closed subalgebra (4.4) of ${\cal A}\Gamma$. Let $\chi_{\langle x\rangle}$ be the characteristic function of a conjugacy class $\langle x\rangle$ of $\Gamma$. Then the pointwise multiplication with $\chi_{\langle x\rangle}$ gives rise to a bounded linear operator 
$$
\begin{array}{cccc}
m_{\langle x\rangle}: & {\cal A}_\omega\Gamma & \to & \ell^1(\Gamma)
\end{array}
\eqno(6.3)
$$
\end{cor}
\begin{proof}
Let $\alpha\in{\cal A}_n\Gamma$ (see 4.1). Then, in the notations of 5.4, 
$$
\parallel m_{\langle x\rangle}(\alpha)\parallel_{\ell^1(\Gamma)}\,=\,\langle\Phi(\vert\alpha\vert),\xi_{h(\langle x\rangle)}\rangle
$$
$$
\leq\,C(n,h(\langle x\rangle))\parallel\vert\alpha\vert\parallel_n\,=\,C(n,h(\langle x\rangle))\parallel\alpha\parallel_n.
$$
\end{proof}

\begin{remark}
In \cite{Lo1} John Lott introduces various delocalized $L^2$-invariants for compact Riemannian manifolds. These are associated to a given conjugacy class of the fundamental group. 
He addresses the question under which conditions these invariants are well defined and of topological 
nature, i.e. independent of the choice of the Riemannian metric. Lott gives an affirmative answer for manifolds whose fundamental group $\Gamma$ satisfies the following conditions: there exists a decreasing (or increasing) family of dense Banach subalgebras $({\cal A}_n\Gamma),\,n\in\N,$ of the reduced group $C^*$-algebra $C^*_r(\Gamma)$ such that
\begin{itemize}
\item[$\bullet$] The intersection (union) of the family of subalgebras is dense and  closed under holomorphic functional in $C^*_r(\Gamma)$.
\item[$\bullet$] Any trace on $\C\Gamma$ which is supported on a fixed conjugacy class of $\Gamma$ extends to a bounded trace 
on the family $({\cal A}_n\Gamma)$.
\end{itemize}
He claims that these conditions are satisfied for groups of polynomial growth and word-hyperbolic groups by referring to the Jolissaint algebra (6.2). This algebra satisfies indeed the first condition, but not the second if $\Gamma$ is nonelementary word-hyperbolic, as was pointed out by Thomas Schick. In the sequel Lott published an erratum \cite{Lo2} where he formulates conditions under which his spectral invariants are still defined and metric independent. 

The algebras constructed in the present paper satisfy both of the conditions above and provide thus exactly the kind of input needed for Lott's argument to work. Therefore his results hold as originally stated in \cite{Lo1}. The rather restrictive additional conditions imposed in \cite{Lo2} are unnecessary.
\end{remark}

\subsection{Homogeneous decompositions of cyclic complexes} 

Let $\Gamma$ be an discrete group. It is well known that the cyclic bicomplex \cite{Co1} of the group ring $\C\Gamma$ possesses a homogeneous decomposition \cite{Ni} 
$$
\begin{array}{ccc}
CC(\C\Gamma) & \simeq & \underset{\langle x\rangle}{\bigoplus}\,CC_{\langle x\rangle}(\C\Gamma)
\end{array}
\eqno(6.4)
$$
where $\langle x\rangle$ runs over the set of conjugacy classes of $\Gamma$. In fact, if the cyclic bicomplex is viewed as complex of algebraic differential forms \cite{Co1}, then 
$$
\begin{array}{cc}
CC(\C\Gamma)\,=\,(\Omega(\C\Gamma),\,\partial), &
CC_{\langle x\rangle}(\C\Gamma)\,=\,(\Omega(\C\Gamma)_{\langle x\rangle},\,\partial)
\end{array}
\eqno(6.5)
$$
where
$$
\begin{array}{ccc}
\Omega(\C\Gamma)_{\langle x\rangle} & = & \{\,\text{linear span of}\,\,\,u_{g_0}du_{g_1}\ldots du_{g_n},\,
g_0g_1\ldots g_n\in\langle x\rangle\,\}
\end{array}
\eqno(6.6)
$$ 
If one wants to calculate suitable versions of cyclic homology for completions of group rings then one encounters the basic difficulty that in general there exists no homogeneous decomposition anymore. There are only two noteworthy exceptions. On the one hand the analytic cyclic bicomplex of the group Banach algebra $\ell^1(\Gamma)$ decomposes as topological direct sum (in the $\ell^1$-sense) of its homogeneous components. On the other hand the closure of the subcomplex $CC_{\langle e\rangle}(\C\Gamma)$ 
corresponding to the conjugacy class of the unit element defines a direct topological summand of the analytic cyclic bicomplex of any unconditional admissible Fr\'echet algebra ${\cal A}\Gamma$ over $\C\Gamma$. This was observed by Connes and Moscovici in \cite{CM} and plays a crucial role in their proof of the Novikov higher signature conjecture for word-hyperbolic groups. It is actually this (partial) homogeneous decomposition which makes the local cyclic cohomology \cite{Pu2} of $\ell^1(\Gamma)$ (and its homogeneous part in the case of ${\cal A}\Gamma$) accessible to calculation 
for some classes of groups (see \cite{Pu3}, \cite{Pu4}). 

Here we establish a weak homogeneous decomposition of the local cyclic bicomplex of sufficiently large unconditional admissible Fr\'echet algebras over the group ring of a word hyperbolic group. 

Recall that an unconditional admissible Fr\'echet algebra over $\C\Gamma$ is called {\bf sufficiently large} \cite{Pu3}, 4.3, if it contains for any $\lambda>1$ the Banach algebra 
$$
\ell^1_\lambda(\Gamma)\,=\,
\{\underset{g}{\sum} a_g\,u_g,\,\underset{g}{\sum}\,\lambda^{l_S(g)}\,\vert a_g\vert<\infty\}.
\eqno(6.7)
$$ 

The group Banach algebra $\ell^1(\Gamma)$ and the Jolissaint algebra ${\mathfrak A}\Gamma$ of a hyperbolic group are "sufficiently large".
 
\begin{prop}
Let $\Gamma$ be a word-hyperbolic group. Let ${\cal A}\Gamma$ be a sufficiently large unconditional admissible Fr\'echet algebra over $\C\Gamma$ and let ${\cal A}_n\Gamma,\,n\in\N,$ be the associated family of minimal Banach subalgebras (see 4.1). For a conjugacy class $\langle x\rangle$ of $\Gamma$ denote by ${\cal CC}^\epsilon_{\langle x\rangle}({\cal A}_n\Gamma)$ the closure 
of the subcomplex $CC_{\langle x\rangle}(\C\Gamma)$  in the entire cyclic bicomplex ${\cal CC}^\epsilon({\cal A}_n\Gamma)$ \cite{Co2}. Then 
$$
\begin{array}{ccc}
"\underset{\underset{n}{\rightarrow}}{\lim}"\,{\cal CC}^\epsilon({\cal A}_n\Gamma) & \simeq & 
"\underset{\underset{\lambda}{\rightarrow}}{\lim}"\,{\cal CC}^\epsilon_{\langle x\rangle}(\ell^1_\lambda\Gamma)
\\
\end{array}
\eqno(6.8)
$$
and
$$
\begin{array}{ccccccccc}
"\underset{\underset{n}{\rightarrow}}{\lim}"\,{\cal CC}^\epsilon({\cal A}_n\Gamma) & \simeq & 
"\underset{\underset{n}{\rightarrow}}{\lim}"\,{\cal CC}^\epsilon_{\langle x\rangle}({\cal A}_n\Gamma)
& \oplus & ? &
\simeq & 
"\underset{\underset{\lambda}{\rightarrow}}{\lim}"\,{\cal CC}^\epsilon_{\langle x\rangle}(\ell^1_\lambda(\Gamma))
& \oplus & ? 
\\
\end{array}
\eqno(6.9)
$$
as ind-complex of topological vector spaces. (For the notion of formal inductive limit see \cite{Pu1}).
\end{prop}

\begin{proof}
Using the ideas of \cite{CM}, this is an immediate consequence of 6.3. 

Let $S$ be a finite symmetric set of generators for $\Gamma$ and let $\parallel-\parallel$ be an unconditional seminorm on ${\cal A}\Gamma$ such that $U=\{x\in{\cal A}\Gamma,\,
\parallel x\parallel<1\}$ is an open unit ball for ${\cal A}\Gamma$. Clearly $U$ is stable under $a\mapsto\vert a\vert$.

By assumption $\ell^1_\lambda(\Gamma)\subset{\cal A}\Gamma$ for every $\lambda>1$, and the closed graph theorem shows that the inclusion is bounded. Thus there exists a constant $C_{21}(\lambda)$ such that 
$$
\parallel u_g\parallel\,\leq\,C_{21}(\lambda)\parallel u_g\parallel_{\ell^1_\lambda(\Gamma)}\,=\,C_{21}(\lambda)\,\lambda^{l_S(g)}
\eqno(6.10)
$$
for all $g\in\Gamma$, and
$$
\begin{array}{ccc}
\underset{l_S(g)\to\infty}{\overline{\lim}}\,\parallel u_g\parallel^{\frac{1}{l_S(g)}} & \leq & 1
\end{array}
\eqno(6.11)
$$
In particular, there exists for each $\lambda>1$ an integer  $N=N(\lambda)$ such that
$\lambda^{-\frac12}\,<\,1-N^{-1}$ and
$$
\begin{array}{ccccc}
\parallel u_g\parallel & \leq & \lambda^{\frac{l_S(g)}{2}} & \text{if} & l_S(g)\geq N.
\end{array}
\eqno(6.12)
$$
Thus, according to the definition of the seminorms 4.1. 
$$
\begin{array}{ccccc}
\parallel u_h\parallel_{N} & \leq & (1-N^{-1})^{-1}\parallel u_h\parallel & \leq  & \lambda^{l_S(h)}
\end{array}
\eqno(6.13)
$$
for elements $h\in\Gamma$ of word length $l_S(h)=N$. Using the submultiplicativity of the seminorm $\parallel-\parallel_N$ one deduces easily for suitable $C_{22}(\lambda)>0$ and all $g\in\Gamma$ the inequality 
$$
\begin{array}{ccc}
\parallel u_g\parallel_{\varphi(n)} & \leq & C_{22}(\lambda)\cdot \lambda^{l_S(g)}
\end{array}
\eqno(6.14)
$$
which means that the identity of $\C\Gamma$ gives rise to a bounded operator 
$$
\begin{array}{ccc}
\ell^1_\lambda(\Gamma) & \to & {\cal A}_N(\Gamma)
\end{array}
\eqno(6.15)
$$
Passing to formal inductive limits one obtains a homomorphism of ind-algebras
$$
\begin{array}{ccc}
"\underset{\underset{\lambda}{\rightarrow}}{\lim}"\,\ell^1_\lambda(\Gamma) & \longrightarrow & 
"\underset{\underset{n}{\rightarrow}}{\lim}"\,{\cal A}_n(\Gamma) 
\end{array}
\eqno(6.16)
$$
which induces corresponding morphisms of the various cyclic chain complexes.

Recall the definition of the norms on the minimal Banach subalgebras ${\cal A}_n\Gamma$ (see 4.1). 
For $\alpha\in\C\Gamma$ one has 
$\parallel\alpha\parallel_n\,=\,\underset{I}{\Inf}(\underset{I}{\sum}\nu_i)$ 
where the infimum is taken over all finite families of real numbers $\nu_i>0$ 
and of elements $\alpha_{ij}\in\C\Gamma$ of norm $\parallel\alpha_{ij}\parallel<1-\frac{1}{n}$ in the linear span of $\{u_g,\,l_S(g)\leq n\}$ satisfying 
$\vert\alpha\vert\,\leq\,\underset{I}{\sum}\nu_i\vert\alpha_{i1}\vert...\vert
\alpha_{i{k_i}}\vert.$

Fix now $n\in\N$ and let $C_{23}(n)>1$ be such that $C_{23}(n)(1-\frac{1}{n})<1-\frac{1}{n+1}$. Then there exists 
some $C_{24}(n)>1$ such that 
$$
\begin{array}{ccc}
\vert\underset{g}{\sum}a_g\,u_g\vert\,\leq\,\underset{I}{\sum}\mu_i\vert\alpha_{i1}\vert...
\vert\alpha_{i{k_i}}\vert & \Rightarrow &
\vert\underset{g}{\sum}a_g\,C_{24}^{l_S(g)}\,u_g\vert\,\leq\,
\underset{I}{\sum}\mu_i\vert C_{23}\alpha_{i1}\vert...\vert C_{23}\alpha_{i{k_i}}\vert
\end{array}
\eqno(6.17)
$$
if the elements $\alpha_{ij}$ lie in the linear span of $\{u_g,\,l_S(g)\leq n\}$. This immediately implies
$$
\begin{array}{ccc}
\parallel\underset{g}{\sum}a_g\,C_{24}(n)^{l_S(g)}\,u_g\parallel_{n+1} & \leq &
\parallel\underset{g}{\sum}a_g\,u_g\parallel_n
\end{array}
\eqno(6.18)
$$
for all elements $\alpha=\underset{g}{\sum}a_g\,u_g\in\C\Gamma$.

After these preliminaries we are ready to prove 6.5. 

Choose $\lambda$ such that $1<\lambda<C_{24}(n).$
Denote by $\pi_{\langle x\rangle}:\Omega \C\Gamma\to\Omega\C\Gamma_{\langle x\rangle}$ the canonical projection. Let $a^i\,=\,\underset{g}{\sum}\,a^i_g\,u_g\in\C\Gamma,\,i=0,...,m.$ Then
$$
\parallel\pi_{\langle x\rangle}(a^0da^1...da^m)\parallel_{\ell^1_\lambda(\Gamma)}\,\leq\,
\underset{g_0,g_1,...,g_m}{\sum}
\vert a^0_{g_0}\vert\vert a^1_{g_1}\vert\ldots\vert a^m_{g_m}\vert\,
\parallel\pi_{\langle x\rangle}(u_{g_0}du_{g_1}...du_{g_m})\parallel_{\ell^1_\lambda(\Gamma)}
$$
$$
=\,\underset{g_0g_1...g_m\in\langle x\rangle}{\sum}
(\vert a^0_{g_0}\vert\lambda^{l_S(g_0)})\,(\vert a^1_{g_1}\vert\lambda^{l_S(g_1)})\ldots(\vert a^m_{g_m}\vert\,\lambda^{l_S(g_m)})
$$
$$
=\,\parallel m_{\langle x\rangle}(\beta_0\beta_1\ldots\beta_m)\parallel_{\ell^1(\Gamma)}
$$
where 
$$
\begin{array}{cccc}
\beta_i & = & \underset{g}{\sum}\,\vert a^i_{g}\vert\parallel\,\lambda^{l_S(g)}\,u_g\,\in\,\C\Gamma, & i=0,...,m.
\end{array}
\eqno(6.19)
$$
and $m_{\langle x\rangle}$ is the canonical projection of $\C\Gamma$ onto $\C\Gamma_{\langle x\rangle}$
According to Corollary 6.3
$$
\parallel\pi_{\langle x\rangle}(a^0da^1...da^m)\parallel_{\ell^1_\lambda(\Gamma)}\,\leq\,
\parallel m_{\langle x\rangle}(\beta_0\beta_1\ldots\beta_m)\parallel_{\ell^1(\Gamma)}
$$
$$
\leq\,C_{25}(n)\parallel\beta_0\beta_1\ldots\beta_m\parallel_{n+1}
$$
$$
\leq\,
C_{25}(n)\parallel\beta_0\parallel_{n+1}\parallel\beta_1\parallel_{n+1}\ldots\parallel\beta_m\parallel_{n+1}
\eqno(6.20)
$$
From (6.18), (6.19) and the choice of $\lambda$ we derive further
$$
\begin{array}{ccc}
\parallel\beta_i\parallel_{n+1} & \leq & \parallel\alpha_i\parallel_n.
\end{array}
\eqno(6.21)
$$
so that finally
$$
\begin{array}{ccccccc}
\parallel\pi_{\langle x\rangle}(a^0da^1...da^m)\parallel_{\ell^1_\lambda(\Gamma)} & \leq &
C_{25}(n)\parallel a_1\parallel_n\ldots\parallel a_m\parallel_n & = &
C_{25}(n)\parallel a^0da^1...da^m\parallel_n
\end{array}
\eqno(6.22)
$$
which establishes the proposition.
\end{proof}

\begin{theorem}
Let $\Gamma$ be a word-hyperbolic group and let ${\cal A}\Gamma$ be a sufficiently large unconditional admissible Fr\'echet algebra over $\C\Gamma$. Let $\langle x\rangle$  be a conjugacy class of $\Gamma$. The notations of 6.5 are understood. Then there is a direct sum decomposition 
$$
\begin{array}{ccccccccc}
{\cal CC}({\cal A}\Gamma) & \simeq & 
"\underset{\underset{n}{\rightarrow}}{\lim}"\,{\cal CC}^\epsilon_{\langle x\rangle}({\cal A}_n\Gamma)
& \oplus & ?  & \simeq & 
{\cal CC}_{\langle x\rangle}(\ell^1\Gamma) & \oplus & ? 
\\
\end{array}
\eqno(6.23)
$$
of the analytic cyclic bicomplex \cite{Pu2} ${\cal CC}({\cal A}\Gamma)$ in the derived ind-category $ind\,{\cal D}$ \cite{Pu2}.
\end{theorem}

\begin{proof}
This results from 6.5 and the fact that, according to \cite{Pu2}, 6.13, the canonical map
$$
\begin{array}{ccccc} 
"\underset{\underset{n}{\rightarrow}}{\lim}"\,{\cal CC}^\epsilon({\cal A}_n\Gamma)
& \longrightarrow &  {\cal CC}({\cal A}\Gamma)\\
\end{array}
\eqno(6.24)
$$
is an isomorphism in the derived ind-category.
\end{proof}

It should be noted that, similarly to the algebraic case \cite{Ni}, the direct factors $"\underset{\underset{\lambda}{\rightarrow}}{\lim}"\,{\cal CC}_{\langle x\rangle}(\ell^1_\lambda(\Gamma))$ are contractible if $\langle x\rangle$ consists of elements of infinite order \cite{Pu3}. 

A rigorous treatment of cyclic complexes of unconditional Fr\'echet algebras over group rings of hyperbolic groups will be the subject of another paper.

\end{document}